\documentclass[10pt]{article}
\usepackage{mathrsfs}
\usepackage{amsmath}
\usepackage{amssymb}
\usepackage{graphicx}
\def \la{\lambda}

\newtheorem{theorem}{\scshape \mdseries  Theorem}[section]
\newtheorem{lemma}[theorem]{\scshape \mdseries  Lemma}
\newtheorem{coro}[theorem]{\scshape \mdseries  Corollary}
\newtheorem{prop}[theorem]{\scshape \mdseries  Proposition}

\begin{document}

\title{Maximum Estrada Index of Bicyclic Graphs\thanks{
      Supported by National Natural Science Foundation of China (11071002), Program for New Century Excellent
Talents in University, Key Project of Chinese Ministry of Education
(210091), Specialized Research Fund for the Doctoral Program of
Higher Education (20103401110002), Science and Technological Fund of
Anhui Province for Outstanding Youth (10040606Y33), Scientific
Research Fund for Fostering Distinguished Young Scholars of Anhui
University(KJJQ1001), Academic Innovation Team of Anhui University
Project (KJTD001B), Fund for Youth Scientific Research of Anhui
University(KJQN1003). }}
\author{Long Wang, Yi-Zheng Fan\thanks{Corresponding author. E-mail addresses:
fanyz@ahu.edu.cn (Y.-Z. Fan), wanglongxuzhou@126.com (L. Wang),
wangy@ahu.edu.cn (Y. Wang)}, Yi Wang
\\
  {\small \it School of Mathematical Sciences, Anhui University, Hefei 230601, P. R. China}
 }
\date{}
\maketitle

\noindent {\bf Abstract}\ \
Let $G$ be a simple graph of order $n$, let $\lambda_{1}(G),\lambda_{2}(G),\ldots, \lambda_{n}(G)$ be the eigenvalues of the adjacency matrix of $G$.
The Esrada index of $G$ is defined as $EE(G)=\sum_{i=1}^{n}e^{\lambda_{i}(G)}$.
In this paper we determine the unique graph with maximum Estrada index among bicyclic graphs with fixed order.


\noindent {\bf MR Subject Classifications:} 05C50

\section{Introduction}
Let $G$ be a simple graph of order $n$ and let $A(G)$ be its adjacency matrix. The {\it eigenvalues of $G$} are referred to the eigenvalues of $A(G)$,
 denoted by $\lambda_{1}(G)\ge \lambda_{2}(G)\ge \cdots \ge \lambda_{n}(G)$. The {\it Estrada index $EE(G)$} of the graph $G$ is defined as $EE(G)=\sum_{i=1}^{n}e^{\lambda_{i}(G)}$.
The Estrada index was first introduced by Estrada \cite{est00} in 2000.
It was found useful in biochemistry and complex networks, see \cite{est02,est04,estro051,estro052,estrora}.
Recently  the Estrada index has been received a lot of attention in mathematics itself. Many bounds have been established for the
Estrada index in \cite{il,pe,zhang,li,zhou}.
Briefly, for a class $\mathcal S$ of graphs, a  graph $G\in \mathcal S$ is called {\it Estrada maximal} if $EE(G)\geq EE(H)$ for any $H \in \mathcal S$.
The Estrada maximal trees subject to one or more graph parameters have been characterized; see \cite{deng,duzh1, il,zhang}.
The unique Estrada maximal unicyclic graph was also determined in \cite{duzh2}.
So, naturally the next problem is to characterize the Estrada maximal graph among all bicyclic graphs of fixed order.
In this paper, we focus on this problem and determine the unique Estrada maximal graphs among all bicyclic graphs of fixed order.

A {\it bicyclic graph} $G=(V,E)$ is a connected simple graph which satisfies $|E|=|V|+1$.
There are two basic bicyclic graphs: $\infty$-graph and $\theta$-graph.
More concisely, an {\it $\infty$-graph}, denoted by $\infty(p,q,l)$, is obtained from two vertex-disjoint cycles $C_p$ and $C_q$ by
connecting one vertex of $C_p$ and one of $C_q$ with a path $P_l$ of length $l-1$ (in the case of $l = 1$, identifying the above two vertices); and
a {\it $\theta$-graph}, denoted by $\theta(p,q,l)$, is a union of three internally disjoint paths $P_{p+1},P_{q+1},P_{l+1}$ of length $p,q,l$ respectively with common end vertices, where $p, q, l \ge 1$ and at most one of them is $1$.
Observe that any bicyclic graph $G$ is obtained from an $\infty$-graph or a $\theta$-graph $G_0$ (possibly) by attaching trees to some of its vertices.
We call $G_0$ the {\it kernel} of $G$.
%
%


\section{Preliminaries and Lemmas}
Let $M_{k}(G)$ be the $k$-th spectral moment of a graph $G$ of order $n$, i.e., $M_{k}(G)=\sum_{i=1}^{n}\lambda_{i}^{k}(G)$.
It is well known that $M_{k}(G)$ is equal to the number of closed walks of length $k$ in $G$. The
following result reveals the connection between the spectral moments
and Estrada index:
$$EE(G)=\sum_{k=0}^{\infty}\frac{M_{k}(G)}{k!}.$$

For any vertices $u$, $v$ and $w$ (not necessary be distinct) in
$G$, we denote by $M_{k}(G;u,v)$  the number of walks in $G$ with
length $k$ from $u$ to $v$, and by $M_{k}(G;u,v,[w])$  the number of
walks in $G$ with length $k$ from $u$ to $v$ which go through $w$.
Denote by $W_k(G;u,v)$ a walk of length $k$ from $u$ to $v$ in $G$, and by $\mathcal{W}_{k}(G;u,v)$ the set of all such walks.
Clearly $M_{k}(G;u,v)=|\mathcal{W}_{k}(G;u,v)|.$

Let $G_{1}$ and $G_{2}$ be two graphs with $u_{1}, v_{1}\in V(G_{1})$ and $u_{2},v_{2}\in V(G_{2})$. We write $(G_{1};u_{1},v_{1})\preceq (G_{2};u_{2},v_{2})$
if $M_{k}(G_{1};u_{1},v_{1})\leq M_{k}(G_{2};u_{2},v_{2})$ for any positive integer $k$. If, in addition,
$M_{k}(G_{1};u_{1},v_{1})<M_{k}(G_{2};u_{2},v_{2})$ for at least one positive integer $k$, then we write $(G_{1};u_{1},v_{1})\prec (G_{2};u_{2},v_{2})$.
Surely $(G_{1};u_{1},v_{1})= (G_{2};u_{2},v_{2})$ implies $M_{k}(G_{1};u_{1},v_{1})= M_{k}(G_{2};u_{2},v_{2})$ for any positive integer $k$.

\begin{lemma} {\em \cite{dul}} \label{per2} Let $G$ be  a graph containing two vertices $u,v$.
Suppose that $w_{i}\in V(G)$ and $uw_{i}\notin E(G),vw_{i}\notin E(G)$ for $i=1,2,\ldots,k$.
Let $E_{u}=\{uw_{i},i=1,2,\ldots,k\}$ and $E_{v}=\{vw_{i},i=1,2,\ldots,k\}$.
Let $G_{u}=G+E_{u}$ and $G_{v}=G+E_{v}$. If $(G;u,u)\prec (G;v,v)$ and $(G;u,w_{i})\preceq (G;v,w_{i})$ for $1\leq i\leq k$,
 then $EE(G_{u})<EE(G_{v})$.
 \end{lemma}

The {\it coalescence} of two vertex-disjoint connected graphs $G,H$, denoted by $G(u) \circ H(w)$, where $u \in V(G)$ and $w \in V(H)$, is obtained by
identifying the vertex $u$ of $G$ with the vertex $w$ of $H$.
A graph is called {\it nontrival} if it contains at least two vertices.

\begin{lemma}{\em \cite{duzh3}}\label{per-cutv}
Let $G$ be a connected graph containing two vertices $u,v$, and let $H$ be a nontrivial connected graph containing a vertex $w$.
If $(G;u,u) \succ (G;v,v)$, then $EE(G(u)\circ H(w))>EE(G(v)\circ H(w))$.
\end{lemma}

\begin{lemma} {\em \cite{duzh3}} \label{per-cute}
Let $H_1$ be a nontrivial connected graph containing a vertex $w$, and let $H_2$ be a connected graph of order at least $3$  containing an pendant edge $uv$, where $v$ is a pendant vertex.
Then $EE(H_1(w) \circ H_2(u)) > EE(H_1(w) \circ H_2(v))$.
%
\end{lemma}


%
%


\begin{lemma}\label{automor}
Let $H_1$ be a connected graph containing two vertices $u,v$, and let $H_2$ be a connected graph disjoint to $H_1$, which contains a vertex $w$.
Let $H'_2$ be a copy of $H_2$, containing the vertex $w'$ corresponding to $w$ of $H_2$.
Let $G=(H_1(u) \circ H_2(w))(v)\circ H'_2(w')$.
If there exists an automorphism $\sigma$ of $H_1$ such that it interchanges $u$ and $v$, then $(G;u,u) =(G;v,v)$ and $(G;u,t)=(G;v,\sigma(t))$ for any vertex $t$ distinct to $u$.

Furthermore, if letting $\bar{H_1}$ be obtained from $H_1$ by adding some edges incident with $v$ but not $u$, letting $\bar{H'_2}$ be obtained from $H'_2$ by adding some vertices or edges such that the resulting graph is connected, and letting $\bar{G}$ be obtained from $G$ by replacing $H_1$ with $\bar{H_1}$ or $H'_2$ with $\bar{H'_2}$, then  $(\bar{G};u,u)\prec(\bar{G};v,v)$ and $(\bar{G};u,t)\prec(\bar{G};v,\sigma(t))$ for any vertex $t$ distinct to $u$.
\end{lemma}

{\bf Proof.}
Surely $\sigma$ induces an automorphism of $G$, and also induces a 1-1 map from $\mathcal{W}_{k}(G;x,y)$ to $\mathcal{W}_{k}(G;\sigma(x),\sigma(y))$ for any $x,y$ and $k$.
The first assertion follows.

Now we prove the second assertion.
Note that
$$
M_k(\bar{G};u,u)=M_k(\bar{G}-v;u,u)+M_k(\bar{G};u,u,[v]),
M_k(\bar{G};v,v)=M_k(\bar{G}-u;v,v)+M_k(\bar{G};v,v,[u]);$$
and
\begin{align*}
M_k(\bar{G}-v;u,u)&=M_k((\bar{H_1}-v)(u) \circ H_2(w);u,u)\\
&=M_k((H_1-v)(u) \circ H_2(w);u,u)\\
&=M_k((H_1-u)(v) \circ H'_2(w');v,v),
\end{align*}
where the last equality holds as $\sigma$ induces an isomorphism between $(H_1-v)(u) \circ H_2(w)$ and $(H_1-u)(v) \circ H'_2(w')$ and interchanges $u$ and $v$.
However, $$M_k(\bar{G}-u;v,v)=M_k((\bar{H_1}-u)(v) \circ \bar{H'_2}(w');v,v).$$
Since $H_1$ is a proper subgraph of $\bar{H_1}$ or $H'_2$ is a proper subgraph of $\bar{H'_2}$,
we have $$M_k(\bar{G}-v;u,u) \le M_k(\bar{G}-u;v,v)$$ with strict inequality for at least one $k$.

For each walk $W \in \mathcal{W}_k(\bar{G};u,u,[v])$,  write it as $W=W_1W_2$, where $W_1$ is the longest subwalk of $W$ from $u$ to $v$, and $W_2$ is the remaining section from $v$ to $u$.  Define a map $f:
\mathcal{W}_k(\bar{G};u,u,[v]) \to \mathcal{W}_k(\bar{G};v,v,[u])$ by $f(W)=W_2W_1$.
One can verify $f$ is an injection, and hence $M_k(\bar{G};u,u,[v])\le M_k(\bar{G};v,v,[u])$.
So we proved $(\bar{G};u,u)\prec(\bar{G};v,v)$.
The proof of $(\bar{G};u,t)\prec(\bar{G};v,\sigma(t))$ can be argued in a similar way. \hfill $\square$

Denote by $N_G(v)$ the set of neighbors of a vertex $v$ in a graph $G$, and by $d_G(v)$ the cardinality of the set $N_G(v)$.

\begin{coro} \label{unicyclic} Let $G$ be a unicyclic graph obtained from a cycle $C$ by attaching some trees on its vertices.
 Assume $u,w$ are two adjacent vertices on the cycle $C$ such that the tree attached at $u$ is a star centered at $v$ with one of its pendant vertices identified with $u$, and the tree attached at $w$ is a star with its center identified with $w$; see Fig. 2.1.
 If $d_{G}(w)\geq d_{G}(v)+1$, then

\noindent{\em (i)} $(G;w,w)\succ (G;v,v)$;

\noindent{\em (ii)} $(G;w,t)\succ (G;v,t)$ for any $t \notin (N_{G}(v) \cup N_G(w) \cup \{w\}) \setminus V(C)$.
\end{coro}

\begin{center}
\vspace{2mm}
  \includegraphics[scale=.5]{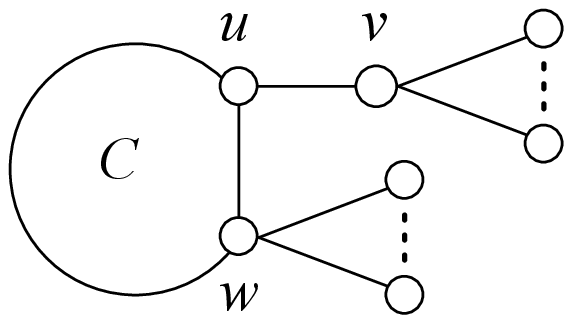}

\vspace{2mm}
{\small Fig. 2.1. An illustration of the graph $G$ in Corollary \ref{unicyclic}}

\end{center}

{\bf Proof.}
Let $G'$ be the graph obtained from $G$ by deleting the edge on the cycle incident to $w$ except $uw$, and deleting $d_G(w)-d_G(v)-1$ pendant vertices of $w$.
Then there exists an automorphism of $G'$ which interchanges $v$ and $w$ together with their pendant vertices, and preserves all other vertices.
Now the assertion follows from the second result of Lemma \ref{automor}. \hfill $\square$

%

\begin{coro} \label{spetheta} Let $G$ be obtained from $\theta(2,2,l)$ by attaching some pendant edges at the vertices of its cycles.
Let $u,v,w,t$ be the vertices as shown in Fig. 2.2.

\noindent{\em (i)} If $d_{G}(w)>2$ and $d_{G}(t)=2$, then $(G;w,w)\succ(G;t,t)$;

\noindent{\em (ii)} If $d_{G}(u)>3$, $d_{G}(v)=3$ and $d_{G}(x)=2$ for any $x\in V(G)\backslash\{u,v,w\}$, then $(G;u,u)\succ(G;v,v)$;

\noindent{\em (iii)} If $d_{G}(u)>3$, $d_{G}(v)=3$ and $d_{G}(x)=2$ for any $x\in V(G)\backslash\{u,v\}$, then $(G;u,u)\succ(G;w,w)$.
\end{coro}

 \begin{center}
 \vspace{2mm}

   \includegraphics[scale=.5]{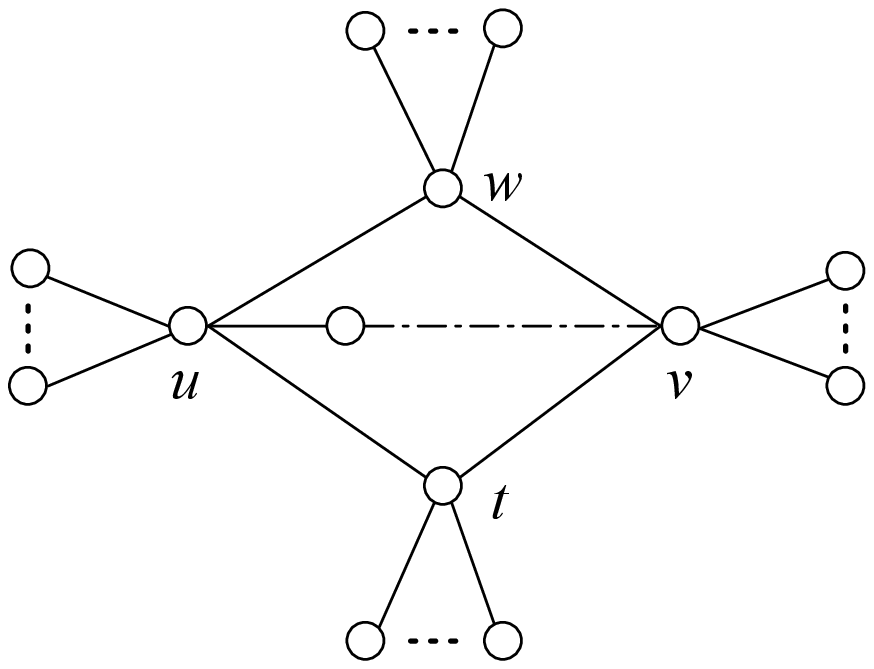}

   \vspace{2mm}

 {\small Fig. 2.2. An illustration of the graph $G$ in Corollary \ref{spetheta}}

\end{center}

{\bf Proof.}
For the assertion (i), let $G'$ be obtained from $G$ by deleting the pendant vertices of $w$.
Then there exists an automorphism $\sigma$ of $G'$ which interchanges $w,t$ and preserves all other vertices.
The assertion follows from Lemma \ref{automor}.
The assertions (ii),(iii) can be argued in a similar way by Lemma \ref{automor}.\hfill$\blacksquare$

\begin{lemma} \label{theta}
Let $G=\theta(p,q,l)$ and let $u,v$ be the two vertices of $G$ with degree $3$ respectively.
Then $(G;u,u)=(G;v,v) \succ (G;w,w)$ for any vertex $w$ distinct to $u$ and $v$.
\end{lemma}
%
%

{\bf Proof.} \
Let $P_{p+1},P_{q+1},P_{l+1}$ be respectively the induced paths of $G$ joining $u$ and $v$.
Define an automorphism $\sigma$ of the graph $G$ as follows:
$\sigma$ interchanges $u$ and $v$, and for each vertex $x$ of the path $P_{p+1}$ (respectively, $P_{q+1},P_{l+1}$), $\sigma(x)$ is also on $P_{p+1}$ (respectively, $P_{q+1},P_{l+1}$) such that the distance between $x$ and $u$ along this path is equal to that between $\sigma(x)$ and $v=\sigma(u)$.
The automorphism $\sigma$ naturally induces a map from $\mathcal W_{k}(G;s,t)$ to $\mathcal W_{k}(G;\sigma(s),\sigma(t))$, such that
$\sigma(W_{k}(G;s,t))=W_{k}(G;\sigma(s),\sigma(t))$ for any $k$ and $s,t\in V(G)$, where each vertex $x$ of $W_{k}(G;s,t)$ is mapped to $\sigma(x)$.
In particular, $\sigma$ is a 1-1 map from $\mathcal W_{k}(G;u,u)$ to $\mathcal W_{k}(G;v,v)$ for any $k$, and hence
$(G;u,u)=(G;v,v)$.

To prove $(G;u,u) \succ (G;w,w)$ for any vertex $w$ distinct to $u$ and $v$,
we only consider the case when $w$ lies on the internal part of the path $P_{l+1}$. The other cases can be proved in a similar way.
Denote $C_{ij}$ the cycle made by $P_{i+1}$ and $P_{j+1}$, where $i \ne j$, and $i,j$ is one of $p,q,l$.
One can easily see that $M_{k}(C_{ql};u,u)=M_{k}(C_{ql};w,w)$.
Thus it suffices to consider the closed walks of length $k$ from $w$ to $w$ that pass though at least one edge of $P_{p+1}$.
Suppose that $W_{k}(G;w,w)$ is such a walk. We decompose the walk $W_{k}(G;w,w)$ into
three parts $W_{1},W_{2},W_{3}$ in a unique way, where
 $W_1$ starts at $w$ and goes along the path $P_{l+1}$ as far as possible, whose terminal point must be $u$ or $v$;
$W_2$ starts at the terminal point of $W_1$, takes the first step and the last step on edges of $C_{pq}$,
and $W-W_2$ contains no edges of $C_{pq}$,  whose terminal point must be $u$ or $v$;
$W_3=W-W_1W_2$.

  Now we construct a map $g$ from $\mathcal{W}_{k}(G;w,w)$ to $\mathcal{W}_{k}(G;u,u)$ in the following
  way. If $W_{2}$ is a $u-u$ walk or $u-v$ walk, $g(W_{1}W_{2}W_{3})=W_{2}W_{3}W_{1}$;
if $W_{2}$ is a $v-u$ walk,  $g(W_{1}W_{2}W_{3})=W_{3}W_{1}W_{2}$;
if $W_{2}$ is a $v-v$ walk, $g(W_{1}W_{2}W_{3})=\sigma(W_{3}W_{1}W_{2})$.
 By directly checking we find that $g$ is an injection. Thus $M_{k}(G;u,u)\geq M_{k}(G;w,w)$ for any $k$. Obviously, $M_{2}(G;u,u)=3>2=M_{2}(G;w,w)$.
 This completes the proof. \hfill $\square$

\section{Main results}
 Denote by $\mathcal G_\infty(n;p,q)$
 the set of all bicyclic graphs of order $n$ which contains an $\infty$-graph as a kernel with two cycles having length $p,q$ respectively.
Denote by $\mathcal G_\theta(n;p,q)$ the set of all bicyclic graphs of order $n$ which contains $\theta(p',q',l')$ as kernel, where $p' \ge q' \ge l'$ and $p'+l'=p, q'+l'=q$.
We first investigate some properties of Estrada maximal graphs in $\mathcal G_\infty(n;p,q)$ or $\mathcal G_\theta(n;p,q)$,
and show that any Estrada maximal graph in $\mathcal G_\infty(n;p,q)$ will have a smaller Estrada index than some graph in $\mathcal G_\theta(n;p,q)$.
Finally we determine the unique  Estrada maximal graph among all bicyclic graphs of fixed order.

\begin{lemma} \label{pend}
If $G$ is an Estrada maximal graph among all bicyclic graphs of order $n$, then $G$ is obtained from its kernel by attaching some pendant edges.
\end{lemma}

{\bf Proof.}
Assume to the contrary, there exists a pendant edge $G$ not attached to its kernel.
Then there is a cut edge $uw$ of $G$ such that $G-uw$ has two components $G_1,G_2$, where $G_1$ contains the vertex $u$ and the kernel of $G$, and $G_2$ is a nontrivial tree containing the vertex $w$.
Removing $G_2$ at $w$ and attaching it to $u$, by Lemma \ref{per-cute} we will arrive at a new bicyclic graph but with larger Estrada index, a contradiction. \hfill $\square$

\begin{theorem} \label{inftymain}
If $G$ is an Estrada maximal graph in $\mathcal G_\infty(n;p,q)$, then $G$ is obtained from $\infty(p,q,1)$ by attaching some pendant edges.
\end{theorem}

{\bf Proof.}
 Suppose $G$ is the Estrada maximal graph in $\mathcal G_\infty(n;p,q)$, and contains $\infty(p,q,l)$ as its kernel.
 By Lemma \ref{pend}, $G$ is obtained from $\infty(p,q,l)$ by attaching some pendant edges.
 We assert $l=1$.
 Otherwise, let $P_l \;(l>1)$ be the path connecting $C_p$ and $C_q$, and let $v_1v_2$ be the starting edge of $P_l$, where $v_1 \in V(C_p)$.
 Write $G=G_1(v_1)\circ G_2(v_1)$,  where $G_1$ contains $C_p$, and $G_2$ contains $C_q$ and the vertex $v_1$ as a pendant vertex.
 Removing $G_1$ at $v_1$ and attaching it to $v_2$, we will arrive at a graph $G' \in \mathcal G_\infty(n;p,q)$.
 However, by Lemma \ref{per-cute},  $EE(G')>EE(G)$, a contradiction.\hfill $\square$

\begin{theorem} \label{thetamain}
If $G$ is an Estrada maximal graph in $\mathcal G_\theta(n;p,q)$, then $G$ is obtained from $\theta(p-1,q-1,1)$ or $\theta(2,2,2)$ by attaching some pendant edges.
\end{theorem}

{\bf Proof.}
Suppose $G$ is an Estrada maximal graph in $\mathcal G_\theta(n;p,q)$, and contains $\theta(p',q',l')$ as its kernel, where $p' \ge q' \ge l'$ and $p'+l'=p, q'+l'=q$.
By Lemma \ref{pend}, $G$ is obtained from $\theta(p',q',l')$ by attaching some pendant edges.
If $l'=1$, or $l'=2$ and $p'=q'=2$, the result follows.
Now assume $l' \ge 2$ and $p' \ge 3$.
Let $u,v,w,t$ be the vertices of $\theta(p',q',l')$ as shown in the left graph of Fig. 3.1.
Without loss of generality, assume $d_G(w) \ge d_G(v)$.
Deleting the edge $tv$ and adding a new edge $tw$, we will arrive at a new graph $\bar{G}$ whose kernel is $\theta(p'-1,q'-1,l')$ as shown in the right graph in Fig. 3.1.
Consider the unicyclic graph $G-tv$.
By Lemma \ref{unicyclic}, $(G-tv;w,w)\succ (G-tv;v,v)$ and $(G-tv;w,t)\succ (G-tv;v,t)$.
So, by Lemma \ref{per2}, $EE(\bar{G})> EE(G)$, a contradiction. \hfill $\square$

%

 \begin{center}
\vspace{2mm}
    \includegraphics[scale=.5]{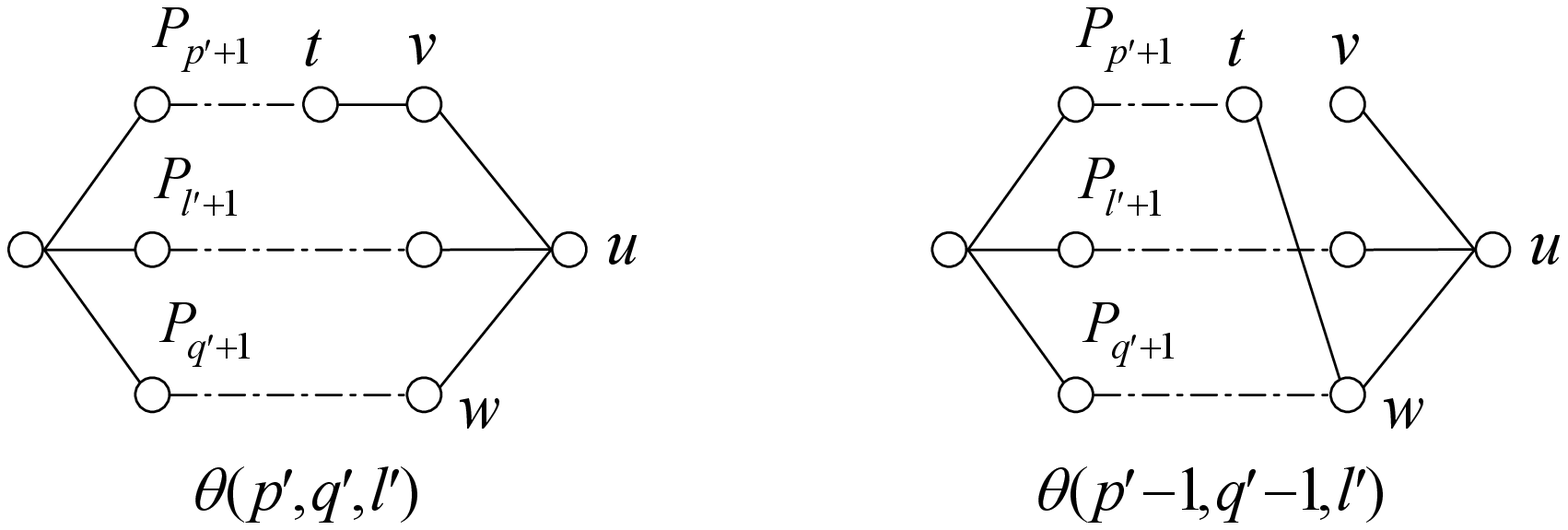}

 \vspace{2mm}
  {\small Fig. 3.1. An illustration of the proof of Theorem \ref{thetamain}}
\end{center}

\begin{lemma} \label{inf-the}
Let $G$ be a bicyclic graph which is obtained from $\infty(p,q,1)$ by attaching some pendant edges to its vertices.
Then their exists a bicyclic graph $\bar{G}$ whose kernel is $\theta(p-1,q-1,1)$ such that $EE(\bar{G}) > EE(G)$.
\end{lemma}

{\bf Proof.}
Let $v,w,t$ be the vertices of $\infty(p,q,1)$ as shown in Fig. 3.2, where $d_G(w) \ge d_G(v)$.
Deleting the edge $tv$ and adding a new edge $tw$, we will arrive at a new graph $\bar{G}$ whose kernel is $\theta(p-1,q-1,1)$ as shown in Fig. 3.2.
Consider the unicyclic graph $G-tv$.
By Lemma \ref{unicyclic}, $(G-tv;w,w)\succ (G-tv;v,v)$ and $(G-tv;w,t)\succ (G-tv;v,t)$.
So, by Lemma \ref{per2}, $EE(\bar{G})> EE(G)$, a contradiction. \hfill $\square$


\begin{center}
\vspace{2mm}

  \includegraphics[scale=.5]{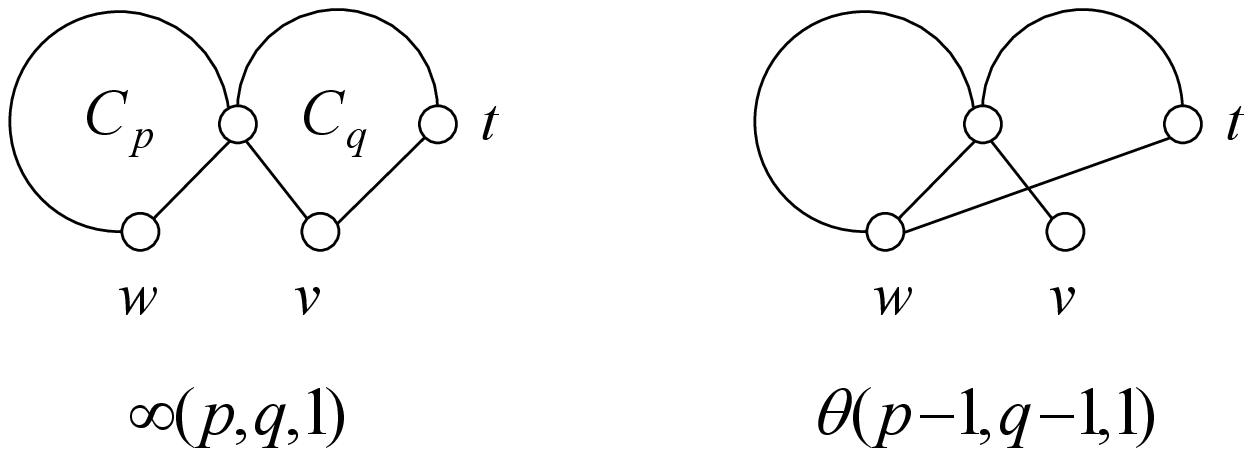}

  \vspace{2mm}
 {\small Fig. 3.2. An illustration of proof of Lemma \ref{inf-the}}
\end{center}

Denote by $\mathbf G_1$ the bicyclic graph of order $n$ obtained from $\theta(2,2,1)$ by attaching $n-4$ pendant edges to one of its vertices of degree $3$,
 and by $\mathbf G_{2}$ the bicyclic graph of order $n$  obtained from $\theta(2,2,2)$ by attaching $n-5$ pendant edges to one of its vertices of degree $3$.

\begin{lemma} \label{max0}
Let $G$ be an Estrada maximal graph among all bicyclic graphs of order $n$.
Then $G$ is either $\mathbf G_1$ or $\mathbf G_{2}$.
\end{lemma}

{\bf Proof.}
By Theorem \ref{inftymain} and Lemma \ref{inf-the}, $G$ must contains a $\theta$-graph as its kernel.
By Lemma \ref{thetamain}, $G$ is obtained from $\theta(p,q,1)$ or $\theta(2,2,2)$ by attaching some pendant edges.
Assume $G$ contains $\theta(p,q,1)$ as its kernel, where $p \ge q$ and $p \ge 3$.
Let $v,w,t$ be the vertices of $\theta(p,q,1)$ as shown in Fig. 3.3, where $d_G(w) \ge d_G(v)$.
Deleting the edge $tv$ and adding a new edge $tw$, we will arrive at a new graph $\bar{G}$ whose kernel is $\theta(p-1,q-1,2)$ as shown in Fig. 3.3.
By a similar discussion to the proof of Lemma \ref{inf-the}, we have $EE(\bar{G})>EE(G)$, a contradiction.
So $G$ is obtained from $\theta(2,2,1)$ or $\theta(2,2,2)$ by attaching some pendant edges.

We next show all the pendant edges of $G$ are attached at a unique vertex $\theta(2,2,1)$ or $\theta(2,2,2)$ with degree $3$,
and hence $G$ is exactly  $\mathbf G_1$ or $\mathbf G_{2}$.
We only prove the case of $G$ having $\theta(2,2,2)$ as the kernel; the other case can be discussed in a similar way.
Let $v_i, i=1,2,\ldots,5$, be the vertices of $\theta(2,2,2)$ as shown in the last graph in Fig. 3.3.
Assume each $v_i$ is attached to $m_i$ pendant edges in the graph $G$, for $i=1,2,\ldots,5$, where $m_i \ge 0$ and $\sum_{i=1}^5=n-5$.

Denote $G=:G(m_1,m_2,m_3,m_4,m_5)$.
If at least two of $m_3,m_4,m_5$ are nonzero, say $m_3>0,m_4>0$, by Corollary \ref{spetheta}(i),
$(G(m_1,m_2,m_3,0,m_5;v_3,v_3) \succ (G(m_1,m_2,m_3,0,m_5;v_4,v_4)$,
and by Lemma \ref{per-cutv} removing all the pendant edges of $G(m_1,m_2,m_3,m_4,m_5)$ at $v_4$ and attaching them to $v_3$,
we will get a graph $G(m_1,m_2,m_3+m_4,0,m_5)$ with a larger Estrada index, a contradiction.
So, at least two of $m_3,m_4,m_5$ are zero, say $m_4=m_5=0$.
Then $G=G(m_1,m_2,m_3,0,0)$.

If both $m_1,m_2$ are nonzero, by Corollary \ref{spetheta}(ii),
$(G(m_1,0,m_3,0,0;v_1,v_1) \succ (G(m_1,0,m_3,0,0;v_2,v_2)$, and by Lemma \ref{per-cutv} removing all the pendant edges of $G(m_1,m_2,m_3,0,0)$ at $v_2$ and attaching them to $v_1$,
we will arrive at graph $G(m_1+m_2,0,m_3,0,0)$ with a larger Estrada index, also a contradiction.
So at least one of $m_1,m_2$ is zero, say $m_2=0$, i.e. $G=G(m_1,0,m_3,0,0)$.
By Lemmas \ref{theta} and \ref{spetheta}, $(G(m_1,0,0,0,0;v_1,v_1) \succ (G(m_1,0,0,0,0;v_3,v_3)$ whether or not $m_1=0$.
If $m_3>0$, by a similar discussion we get $EE(G(m_1+m_3,0,0,0,0))>EE(G(m_1,0,m_3,0,0))$, a contradiction.
So $G=G(m_1,0,0,0,0)$, and the result follows. \hfill $\square$

\begin{center}
\vspace{3mm}

  \includegraphics[scale=.5]{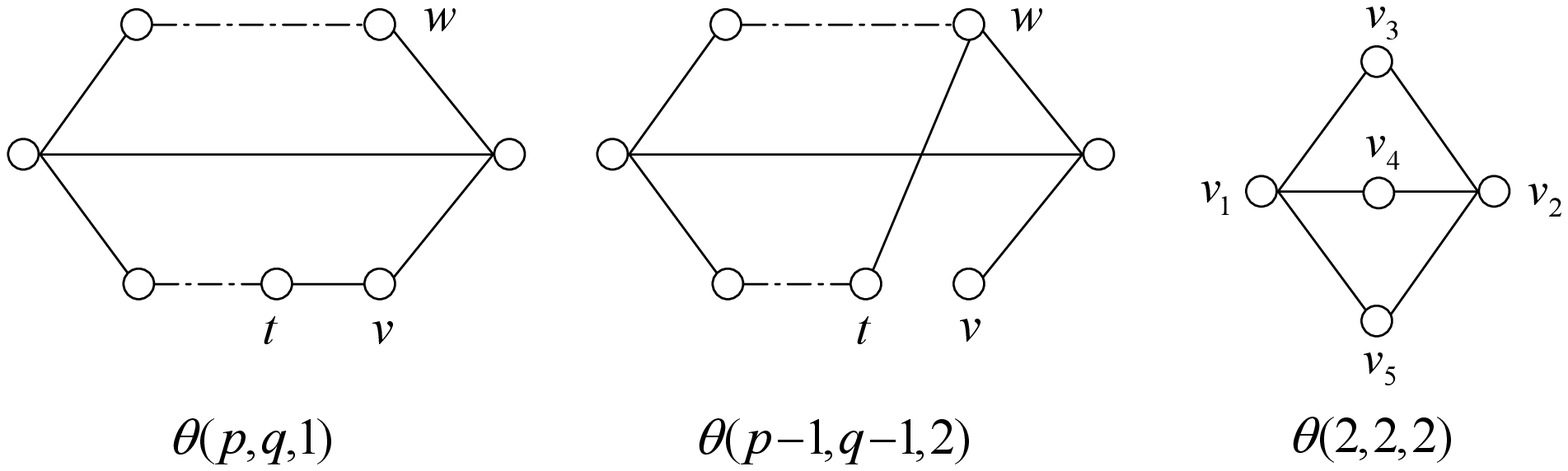}

 \vspace{3mm}
 {\small Fig. 3.3. An illustration of proof of Lemma \ref{max0}}

\end{center}

Finally we determine which is larger between $EE(\mathbf G_1)$ and $EE(\mathbf G_2)$.
Denote by $\phi(G,x)$ the characteristic polynomial of the adjacency matrix of a graph $G$.

 \begin{lemma} {\em \cite{cve}} \label{poly}
 Let $G$ be a graph containing a vertex $v$, and let $\mathcal{C}(v)$ be the set of cycles containing $v$.
Then
 $$\phi(G,x)=x \phi(G-v,x)-\sum_{w \in N_G(v)}\phi(G-v-w,x)-2\sum_{Z\in \mathcal{C}(v)}\phi(G-V(Z),x).$$
\end{lemma}

\begin{prop} \label{G12} $EE(\mathbf G_1)> EE(\mathbf G_2)$ for $n\geq 5$.
\end{prop}

{\bf Proof.}
By Proposition \ref{poly}, $\phi(\mathbf G_{1},x)=x^{n-4}f(x)$, $\phi(\mathbf G_{2},x)=x^{n-4}g(x)$, where
$$f(x)=x^{4}-(n+1)x^{2}-4x+2(n-4), g(x)=x^{4}-(n+1)x^{2}+3(n-5).$$
By a direct calculation, $EE(\mathbf G_{1})>EE(\mathbf G_{2})$ when $5\leq n\leq 22$.
Now assume $n\geq 23$.
Since $f(\sqrt{n-1})=-6-4\sqrt{n-1}<0$, $\la_1(\mathbf G_{1})>\sqrt{n-1}$.
On the other hand, as $g(x)$ is increasing for $x>\sqrt{\frac{n+1}{2}}$,  $g\left(\sqrt{n-\frac{3}{2}}\right)=\frac{n}{2}-\frac{45}{4}>0$ when $n\geq 23$,
 which implies $\la_1(\mathbf G_{2})<\sqrt{n-\frac{3}{2}}$ when $n\geq 23$.

Let $u,v$ be the vertices of $\mathbf G_1$ and $\mathbf G_2$ both with maximal degree, respectively.
 The graph $\mathbf G_{1}-u$ has eigenvalues $\pm\sqrt{2}$ and $0$ with multiplicity $n-3$, and
  the graph $\mathbf G_{2}-v$ has eigenvalues $\pm\sqrt{3}$ and $0$ with multiplicity $n-3$.
  By interlacing property of the eigenvalues of $A( \mathbf G_{1}-u)$ and $A(\mathbf G_{1}) $ (or see \cite{cve}),
  $\la_i(\mathbf G_{1}) \ge \la_i(\mathbf G_{1}-u)$ for $i=2,3,\ldots,n-1$.
  So
$$
 EE(\mathbf G_{1}) =\sum_{i=1}^n e^{\la_i(\mathbf G_{1})}  > e^{\lambda_{1}(\mathbf G_{1})}+\sum_{i=2}^{n-1} e^{\la_i(\mathbf G_{1}-u)}
   > e^{\sqrt{n-1}}+(n-3)+e^{-\sqrt{2}}.
$$
Similarly, by the fact $\la_i( \mathbf G_{2}) \le  \la_{i-1}( \mathbf G_{2}-v)$ for $i=2,3,\ldots,n$,
$$EE(\mathbf  G_{2})\leq e^{\la_1( \mathbf G_{2})}+\sum_{i=2}^{n} e^{\la_i(\mathbf G_{2}-v)}<e^{\sqrt{n-\frac{3}{2}}}+e^{\sqrt{3}}+(n-3)+e^{-\sqrt{3}}.$$
Noting that $e^{\sqrt{n-1}}>e^{\sqrt{n-\frac{3}{2}}}+e^{\sqrt{3}}$ for $n\geq 23$, so we get the result. \hfill $\square$

By Lemma \ref{max0} and Proposition \ref{G12}, we get the main result of this paper.

\begin{theorem} \label{max}
Let $G$ be a bicyclic graph of order $n$.
Then $EE(G) \le EE(\mathbf G_1)$, with equality if and only if $G=\mathbf G_1$.
\end{theorem}


\small
\begin{thebibliography}{90}
\bibitem{cve} D. Cvetkov\'{i}c, M. Doob, H. Sachs, {\it Spectra of Graphs-Theory and Application}, Academic Press, New York, 1980.
\bibitem{deng} H. Deng, A proof of a conjecture on the Estrada index, {\it MATCH Commun. Math. Comput. Chem.}, 62 (2009) 599-606.
\bibitem{dul} Z. Du, Z. Liu, On the Estrada and Laplacian Estrada indices of graphs, {\it Linear Algebra Appl.}, 435 (2011) 2065-2076.
\bibitem{duzh1} Z. Du, B. Zhou, The Estrada index of trees, {\it Linear Algebra Appl.}, 435 (2011) 2462-2467.
\bibitem{duzh3} Z. Du, B. Zhou, On the Estrada index of graphs with given number of cut vertices, {\it Electron. J. Linear Algebra}, 22 (2011) 586-592.
\bibitem{duzh2} Z. Du, B. Zhou, The Estrada index of unicyclic graphs, {\it Linear Algebra Appl.}, 436 (2012) 3149-3159.
\bibitem{est00}  E. Estrada, Characterization of 3D molecular structure, {\it Chem. Phys. Lett.}, 319 (2000) 713-718.
\bibitem{est02} E. Estrada, Characterization of the folding degree of proteins, {\it Bioinformatics}, 18 (2002) 697-704.
\bibitem{est04} E. Estrada, Characterization of the amino acid contribution to the folding degree of proteins, {\it Proteins}, 54 (2004) 727-737.
\bibitem{estro051} E. Estrada, J. A. Rodr\'{i}guez-Val\'{a}zquez, Subgraph centrality in complex networks, {\it Phys. Rev. E.}, 71 (056103) (2005) 1-9.
\bibitem{estro052} E. Estrada, J. A. Rodr\'{i}guez-Val\'{a}zquez, Spectral measures of bipartivity in complex networks, {\it Phys. Rev. E}, 72 (046105) (2005) 1-6.
\bibitem{estrora} E. Estrada, J. A. Rodr\'{i}guez-Val\'{a}zquez, M. Rand\'{i}c, Atomic branching in molecules, {\it Int. J. Quantum Chem.}, 106 (2006) 823-832.
\bibitem{il} A. Il\'{i}c, D. Stevanov\'{i}c, The Estrada index of chemical trees, {\it J. Math. Chem.}, 47 (2010) 305-314.
\bibitem{li} W. Li, A. Chang, On the trees with maximum nullity, {\it MATCH Commun. Math. Comput. Chem.}, 2006, 56(3) 501-508.
\bibitem{pe} J. A. de la Pe\~{n}a, I. Gutman, J. Rada, Estimating the Estrada index, {\it Linear Algebra Appl.}, 427 (2007) 70-76.
\bibitem{zhang} J. Zhang, B. Zhou, J. Li, On Estrada index of trees, {\it Linear Algebra Appl.}, 434 (2011) 215-223.
\bibitem{zhou} B. Zhou, On Estrada index, {\it MATCH Commun. Math. Comput. Chem.}, 60 (2008) 485-492.
\bibitem{zhoutr} B. Zhou, N. Trinajst\'{i}c, Estrada index of bipartite graphs, {\it Int. J. Chem. Model.}, 1 (2008) 387-394.



\end{thebibliography}
\end{document}